\documentclass[11pt]{amsart}
\usepackage{amsmath,amssymb,amscd}
\usepackage{graphicx}
\usepackage{a4wide}

\newcommand{\RR}{\mathbb{R}}
\newcommand{\ZZ}{\mathbb{Z}}
\newcommand{\QQ}{\mathbb{Q}}
\newcommand{\LL}{\mathcal L}
\newcommand{\gL}{\varLambda}
\newcommand{\eins}{\boldsymbol{1}}
\newcommand{\dens}{\mathrm{dens}}
\newcommand{\vol}{\mathrm{vol}}
\newcommand{\oplam}{\mbox{\Large $\curlywedge$}}

\theoremstyle{definition}
\newtheorem{thm}{Theorem}

\newtheorem{coro}{Corollary}

\begin{document}

\title[Homometric model sets]{Homometric model sets and window covariograms}

\author{Michael Baake}
\address{Fakult\"at f\"ur Mathematik, Universit\"at Bielefeld,
Postfach 100131, 33501 Bielefeld,  Germany}
\email{\texttt{mbaake@math.uni-bielefeld.de}}
\urladdr{\texttt{http://www.math.uni-bielefeld.de/baake}}
\author{Uwe Grimm}
\address{Dept.\ of Mathematics, The Open University, 
Walton Hall, Milton Keynes MK7 6AA, UK}
\email{u.g.grimm@open.ac.uk}
\urladdr{http://mcs.open.ac.uk/ugg2/}

\begin{abstract}
Two Delone sets are called homometric when they share the same
autocorrelation or Patterson measure.  A model set $\varLambda$ within
a given cut and project scheme is a Delone set that is defined through
a window $W$ in internal space. The autocorrelation measure of
$\varLambda$ is a pure point measure whose coefficients can be
calculated via the so-called covariogram of $W$. Two windows with the
same covariogram thus result in homometric model sets. On the other
hand, the inverse problem of determining $\varLambda$ from its
diffraction image ultimately amounts to reconstructing $W$ from its
covariogram. This is also known as Matheron's covariogram problem. It
is well studied in convex geometry, where certain uniqueness results
have been obtained in recent years. However, for non-convex windows,
uniqueness fails in a relevant way, so that interesting applications
to the homometry problem emerge.  We discuss this in a simple setting
and show a planar example of distinct homometric model sets.
\end{abstract}

\maketitle

\section{Introduction}

Quasicrystals form an interesting class of solids whose structure
determination requires a substantially extended setting in comparison
to classical crystallography, see \cite{Janot,Steurer} and references
therein. In particular, one first has to extract from the diffraction
data the higher-dimensional embedding space and an appropriate lattice
in it, followed by the reconstruction of the window (or acceptance
domain) for the description of the atomic positions.

In what follows, we begin to analyze the structure of this inverse
problem from a mathematical perspective. To this end, we assume a
perfect diffraction experiment and the appropriateness of 
describing the underlying structure as a model set. The first
step, finding the Fourier module and then the resulting cut and
project scheme, seems relatively straight-forward. We thus assume
this step already done, see also \cite{BL} and further remarks
below for more.

Consequently, we need to determine the window to complete the picture.
Here, we discuss the question how and to what extent it is specified
by the diffraction image. Since all principal aspects show up already
for a mono-atomic system, we make this simplifying assumption, too.
We are well aware of the fact that real quasicrystals require the
reconstruction of several windows, but this adds a technical
complication, not a principal hurdle, to the picture described below.
Complementary reconstruction methods on the basis of discrete
tomography \cite{G} look also promising in the context of quasicrystals
\cite{H+}, but are not considered here.

A Delone set $\gL$ is described by the corresponding Dirac comb
$\delta_{\gL}:=\sum_{x\in\gL}\delta_x$, where $\delta_x$ is the
normalized point (or Dirac) measure at $x$, see \cite{B} for details.
A perfect diffraction image of $\gL$, as described by the positive
measure $\widehat{\gamma}^{}_{\gL}$, uniquely determines its inverse
Fourier transform, which is the autocorrelation (or Patterson) measure
$\gamma^{}_{\gL}$.  Our starting point is thus the (hypothetically
complete) knowledge of $\gamma^{}_{\gL}$, where we shall work with
infinite point sets $\gL$ to avoid the extra layer of complication
from finite apertures. The remaining task is then to determine a
window $W$ from this information, which leads to Matheron's
covariogram problem \cite{Bi2}.

This is an example of a class of inverse problems where one aims at
the reconstruction of a finite or compact set in Euclidean space from
its (possibly weighted) difference set. As originally observed by
Patterson \cite{Patt} for finite point sets in $\RR$, there need not
be a unique solution to the reconstruction problem. To capture this
ambiguity, two point sets are called \emph{homometric} when they share
the same (weighted) difference set, see \cite{Z1} and references
therein for further examples. 

Similar questions emerge for non-empty compact subsets
$K\subset\RR^{d}$, where the \emph{covariogram}
$g^{}_{K}(x)=\vol\bigl(K\cap (x+K)\bigr)$ encapsulates the difference
information, and one tries to determine $K$ from the knowledge of
$g^{}_{K}$, see \cite{Bi2,GGZ}. Furthermore, one is also interested in
similar concepts for \emph{infinite} point sets in general, and Delone
sets in particular. Here, two Delone sets are called \emph{homometric}
when they share the same autocorrelation measure, to be defined in
detail below. There is an interesting connection between the homometry
of Delone sets and the covariogram problem for the class of model sets
(also known as cut and project sets) \cite{M1}, which is explored
below.

The paper is organised as follows. The covariogram and its elementary
properties are introduced in the next section, while
Section~\ref{auto-sec} recalls the basic facts about the
autocorrelation of a regular model set and links it with the
covariogram. Within a given cut and project scheme, the homometry
problem thus becomes equivalent to Matheron's covariogram problem
\cite{Bi2}.  Section~\ref{example} shows a planar example of two
different sets with the same covariogram, which give rise to two
distinct homometric model sets.

\section{Properties of the covariogram}

Let $K\subset\RR^d$ be a non-empty compact set which is
the closure of its interior, $K=\overline{K^{\circ}}$. The function
\begin{equation} \label{covario-1}
   g^{}_{K} (x) \; := \; \vol\bigl(K\cap(x+K)\bigr)
\end{equation}
is called the \emph{covariogram} of the set $K$.  This defines a
real-valued function on $\RR^d$ which is inversion symmetric,
$g^{}_{K} (-x) = g^{}_{K} (x)$, and satisfies
\[
\sup_{x\in\RR^d}\, g^{}_{K} (x) \; = \;  g^{}_{K} (0) \; = \; \vol(K).
\]
The function $g^{}_{K} (x)$ is continuous on $\RR^d$, see
\cite[Thm.~3.1]{CJ}.  Moreover, $g^{}_{K}$ is a positive definite
function with compact support. The latter is the so-called
\emph{difference body}
\begin{equation} \label{covario-2}
   \mathrm{supp} (g^{}_{K}) \; = \; K - K 
   \; := \; \{ x-y \mid x,y \in K\} \, ,
\end{equation}
which is always inversion symmetric. In fact, the entire
covariogram of $K$ equals that of $-K$ as well as that of
any translate $K+t$. This means that $g^{}_{K}$ can determine
$K$ at best up to translations and inversion. As we shall
see below, even this is generally not the case.

If $\eins^{}_{K}$ denotes the characteristic function
of $K$, the function  $g^{}_{K} (x)$ is given by the convolution
\begin{equation} \label{covario-3}
   g^{}_{K} (x) \; = \; 
   \bigl(\eins^{}_{K} * \eins^{}_{-K}\bigr) (x) \, .
\end{equation}
The Fourier transform thus satisfies
\begin{equation} \label{covario-4}
   \widehat{g}^{}_{K} (k) \; = \;
   \big\lvert \widehat{\eins}^{}_{K} (k) \big\rvert^2 ,
\end{equation}
which is a positive, analytic function that vanishes at $\infty$.  This
relation is the reason why, if $K$ is itself inversion symmetric in
the sense that $-K=t+K$ for a suitable translation $t$, $\eins^{}_{K}$ can
be reconstructed from the knowledge of $g^{}_{K}$, up to translation
and inversion \cite{CJ}.

\begin{figure}
\includegraphics[width=0.7\textwidth]{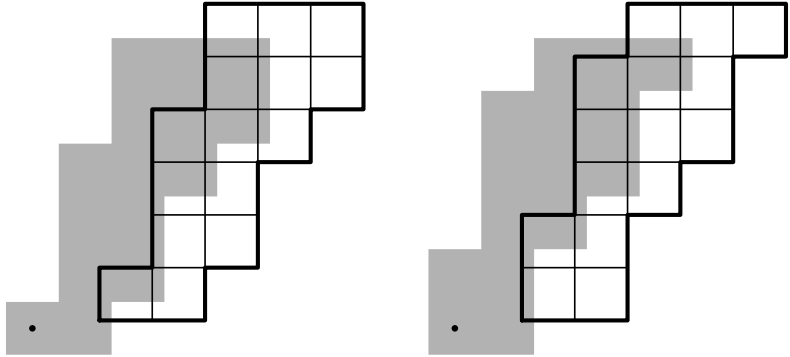}
\caption{Example of two homometric polyominoes, each with an equally
shifted copy overlayed. This illustrates that and how the two
intersection areas coincide. For their use as windows in the cut and
project scheme, we choose the marked dot as the origin in internal
space.\label{fig1}}
\end{figure}

If $K$ is a convex polytope in dimension $d\le 3$, it is determined by
$g_{K}^{}$, in the sense mentioned above, see \cite{Bi1,Bi2,Bi3} and
references therein. Other examples include convex sets in the plane
with $C^{2}$-smooth boundaries, and large other classes (e.g., non-$C^1$
bodies and bodies that are not strictly convex), though the general
claim is still open, see \cite{Bi2}.
However, starting with dimension $4$, uniqueness fails even within the
class of convex polytopes, compare \cite{Bi2,GGZ}.  In general, the
reconstruction of $K$ from the knowledge of $g^{}_{K}$ is a difficult
problem, which is beyond our scope here. However, an interesting
example of two polyominoes with the same covariogram was constructed
in \cite{GGZ} and is shown in Figure~\ref{fig1}.

\section{Homometry of model sets} \label{auto-sec}

In this section, we use the setting of Euclidean spaces for
simplicity. This covers the most common situation in practice; we
refer the reader to \cite{M1,Martin2} for the more general setting.

The basic cut and project scheme \cite{M1} requires a number of
spaces, sets and mappings as follows.
\begin{equation}\label{candp}
\begin{CD}
\RR^{n} @<\pi_{}^{}<< \RR^{n} \!\times\! \RR^{m} 
  @>\pi_{\text{int}}^{}>> \RR^{m} \\
\bigcup @. \bigcup\makebox[0pt][l]{\,\text{\footnotesize lattice}} 
@. \bigcup\makebox[0pt][l]{\,\text{\footnotesize dense}}\\
L @<{1-1}<< \LL @>>> L^{\star} \\
|| @. @. || \\
L @. \hspace*{-9ex}\xrightarrow{\hspace*{11.5ex}\star\hspace*{11.5ex}}
\hspace*{-9ex} @. L^{\star}
\end{CD}\medskip
\end{equation}
Here, $\RR^{n}$ and $\RR^{m}$ are called the direct (or physical)
space and the internal space, respectively. The mappings $\pi$ and
$\pi_{\text{int}}^{}$ denote the canonical projections. The set $\LL$
is a point lattice in $\RR^{n}\!\times\!\RR^{m}=\RR^{n+m}$ such that
$\pi$ is $1-1$ between $\LL$ and $L=\pi(\LL)$ and that
$L^{\star}=\pi_{\text{int}}^{}(\LL)$ is dense in $\RR^{m}$. Finally,
the $\star$-map is defined on $L$ via
$t^{\star}=\pi_{\text{int}}^{}({(\pi |_{\LL}^{})}^{-1}(t))$ for all
$t\in L$.  It has a unique extension to the rational span $\QQ L$ of
$L$, but not to all of $\RR^n$.

To define a (regular) model set on the basis of the cut and project
scheme \eqref{candp}, we need to specify a window $W\subset\RR^m$.  We
assume that $W$ is a compact set with non-empty interior and the
property that $W=\overline{W^\circ}$. Regularity refers to the extra
requirement that the boundary $\partial W$ has measure $0$ in internal
space.  In this setting, we obtain a (regular) model set $\gL$ as
\begin{equation} \label{model-set-1}
   \gL \; = \; t + \oplam(W) \, ,
\end{equation}
with any $t\in\RR^n$ and $\oplam(W) := \{ x\in L \mid x^{\star} \in W \}$.

Let us now consider the autocorrelation $\gamma^{}_{\gL}$ of the Dirac
comb $\delta_\gL = \sum_{x\in\gL} \delta_x$, which is defined
as
\begin{equation} \label{auto-1}
   \gamma^{}_{\gL} \; := \;
   \lim_{r\to\infty} \frac{1}{\vol (B_r)}\,
   \delta_{\gL\cap B_r} * \delta_{-\gL\cap B_r} ,
\end{equation}
where $B_r$ denotes the open ball of radius $r$ around $0$.  The limit
exists by general properties of model sets, compare \cite{Martin2}.
Moreover, it is independent of the translation $t$ in
\eqref{model-set-1}, wherefore we set $t=0$ from now on without loss
of generality. It is well known that the autocorrelation has the
explicit form
\begin{equation} \label{auto-2}
   \gamma^{}_{\gL} \; := \;
   \sum_{x\in\gL-\gL} \eta(x)\, \delta_x
\end{equation}
where $\gL-\gL$ is a uniformly discrete (and
hence countable) subset of $\RR^n$ and the autocorrelation
coefficient $\eta(x)$ is given by \cite[Prop.\ 3]{M2}
\begin{equation} \label{auto-3}
\begin{split}
   \eta(x) &\; = \; \dens(\gL)\,
   \frac{\vol \bigl(W\cap(W-x^{\star})\bigr)}{\vol (W)} \\[1mm]
   &\; = \; \dens(\LL)\, \vol \bigl(W\cap(W-x^{\star})\bigr) .
\end{split}
\end{equation}
This expresses the autocorrelation of $\gL$ in terms of the
covariogram of the window $W$, see \cite{BG,DCG} for related problems
that can be expressed this way.

Let us recall that two Delone sets are called \emph{homometric} when
they share the same autocorrelation, see \cite{Z2} for an example in
the setting of quasicrystals. Note that this definition, due to the
infinite volume limit in Eq.~\eqref{auto-1}, disregards contributions
from point sets of density $0$. For simplicity, we restrict ourselves
to Delone sets with unique autocorrelations here (meaning independence
from the selected averaging van Hove sequence, compare \cite{Martin2}).
In particular, this is the situation for all model sets.

\begin{thm} \label{main-thm}
   Two model sets obtained from the same cut and project
   scheme are homometric if and only if the defining windows
   share the same covariogram.
\end{thm}

\begin{proof}  
If the two model sets, $\gL$ and $\gL'$, possess the same
autocorrelation, the coefficient $\eta(x)$ is valid for both of
them. Since the cut and project scheme is fixed, the factor
$\dens(\LL)$ in Eq.~\eqref{auto-3} is the same for $\gL$ and
$\gL'$. Consequently, the covariograms of the windows coincide on all
points $x\in\gL-\gL$. By the general setting of a cut and project
scheme, $\gL^{\star}-\gL^{\star}$ is a dense subset of $W-W$, so that
the covariogram, which is a continuous function \cite{CJ}, is completely
determined on $W-W$, while it vanishes on the complement. As the same
argument applies to $W'-W'$, the covariograms of $W$ and $W'$
coincide.

The converse claim is obvious from Eq.~\eqref{auto-3}. 
\end{proof}

It is a well established result, see \cite{Martin2} and references
therein, that the diffraction measure of a regular model set $\gL$ is
a pure point measure, meaning that it is a countable sum of weighted
Bragg (or Dirac) peaks. It is given by $\widehat{\gamma}_{\gL}^{}$,
the Fourier transform of the autocorrelation measure
$\gamma_{\gL}^{}$. Since both $\gamma_{\gL}^{}$ and
$\widehat{\gamma}_{\gL}^{}$ are tempered measures, and since Fourier
transform is invertible on this class of measures, the following
consequence is immediate.

\begin{coro}
Homometric model sets from the same cut and project scheme have the
same diffraction measure. Conversely, kinematic diffraction cannot
discriminate between homometric model sets. \qed
\end{coro}

So far, the formulation was based upon a fixed cut and project scheme.
In general, however, there are different possibilities to generate a
given model set. This also implies that two regular model sets
obtained from \emph{different} cut and project schemes can still be
homometric. If so, they have the same autocorrelation $\gamma$ by
definition, and thus the same diffraction measure $\widehat{\gamma}$.
Since regular model sets are pure point diffractive \cite{H,Martin2},
$\widehat{\gamma}$ comprises a countable sum of Dirac measures with
positive intensity. The $\ZZ$-span of their positions is the so-called
\emph{Fourier module}. Via minimal embedding, it defines an
essentially unique cut and project scheme, whose dual in the sense of
\cite{M1} is a natural setting to describe both model sets in the same
cut and project scheme, see \cite{BL} for details. In this sense, the
restriction in Theorem~\ref{main-thm} is not essential.

\section{A simple planar example} \label{example}

\begin{figure}
\includegraphics[width=0.6\textwidth]{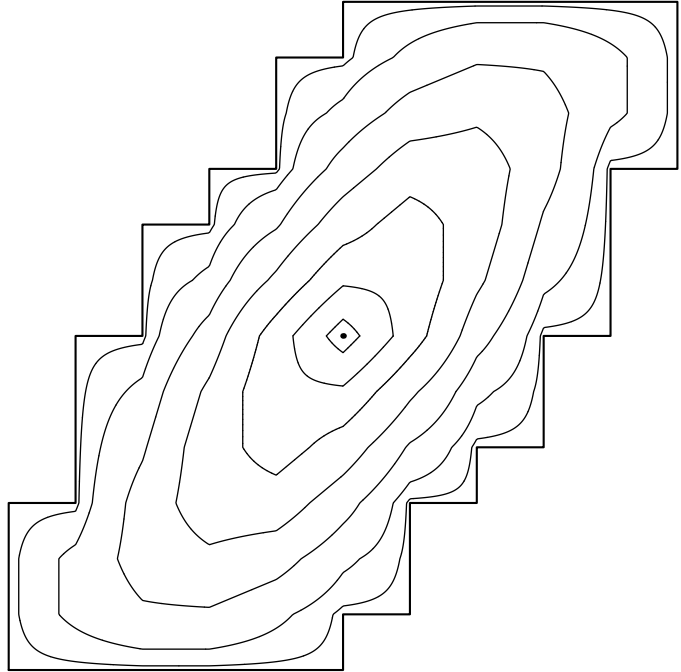}
\caption{Contour plot of the covariogram of both $P_1$ and $P_2$.  It
is inversion symmetric with respect to the origin, the latter marked
by a dot.  The outer polygonal line is the boundary of the difference
body, which is the support of the covariogram, see
Eq.~\eqref{covario-2}.
\label{fig2}}
\end{figure}

A simple example of two distinct homometric polygons was given in
\cite{GGZ}. They are shown in Figure~\ref{fig1}, together with some
typical intersection pattern.  The polygons $P_1$ and $P_2$ are
polyominoes, which both consist of $15$ elementary squares. Their
joint covariogram, supported on the inversion symmetric difference
body, is sketched in Figure~\ref{fig2} as a contour plot.

To turn this pair of polyominoes into a pair of homometric model sets,
we use them as windows for $L=\ZZ[\xi]\subset\RR^2$ with
$\xi=\exp(2\pi i/8)$, with the internal space $\RR^2$ and a $\star$-map
defined by a suitable algebraic conjugation \cite{P}. Here, we use the
one induced by $\xi\mapsto\xi^5$.  With 
\[
   L \; = \; \langle 1,\xi\rangle^{}_{\ZZ[\sqrt{2}\,]} \, ,
\] 
this choice amounts to the replacement $\sqrt{2}\mapsto-\sqrt{2}$ on
the level of the coordinates.

The model set formulation is thus based on a lattice $\LL\subset\RR^4$.
Its basis matrix reads
\begin{equation}\label{basismatrix}
B_{\LL}^{} \; = \; \left(
\begin{array}{rrrr} 
1 &  1/\sqrt{2} & 0 & -1/\sqrt{2}  \\
0 &  1/\sqrt{2} & 1 &  1/\sqrt{2}  \\[0.5ex] 
1 & -1/\sqrt{2} & 0 &  1/\sqrt{2}  \\
0 & -1/\sqrt{2} & 1 & -1/\sqrt{2}
\end{array}\right),
\end{equation}
where the columns contain the coordinates of the basis vectors. The
first two rows show the coordinates in direct space, the remaining
rows those in internal space.  Since $\lvert \det(B_{\LL}^{}) \rvert =
4$, the lattice $\LL$ has density $\dens(\LL) = 1/4$. In this setting,
we define the two model sets
\begin{equation}
\gL_{1} \; = \; \oplam(P_1)  \quad\mbox{and}\quad
\gL_{2} \; = \; \oplam(P_2) .
\end{equation}
Both windows are placed in internal space such that their lower left
corner has coordinates $(-1/2,-1/2)$. This way, the corresponding
model sets are generic, meaning that $L^{\star}\cap\partial P_j
=\varnothing$ for $j\in\{1,2\}$.  They have density
\[
\dens(\gL_1)\; =\; \dens(\gL_2) \; =\; \tfrac{15}{4}
\]
and are homometric by construction, due to Theorem~\ref{main-thm}.
Note that $\gL_1$ and $\gL_2$ are \emph{not}\/ in the same LI class,
compare \cite{B}, and that they differ in points of positive density,
see below.

\begin{figure}
\includegraphics[width=0.65\textwidth]{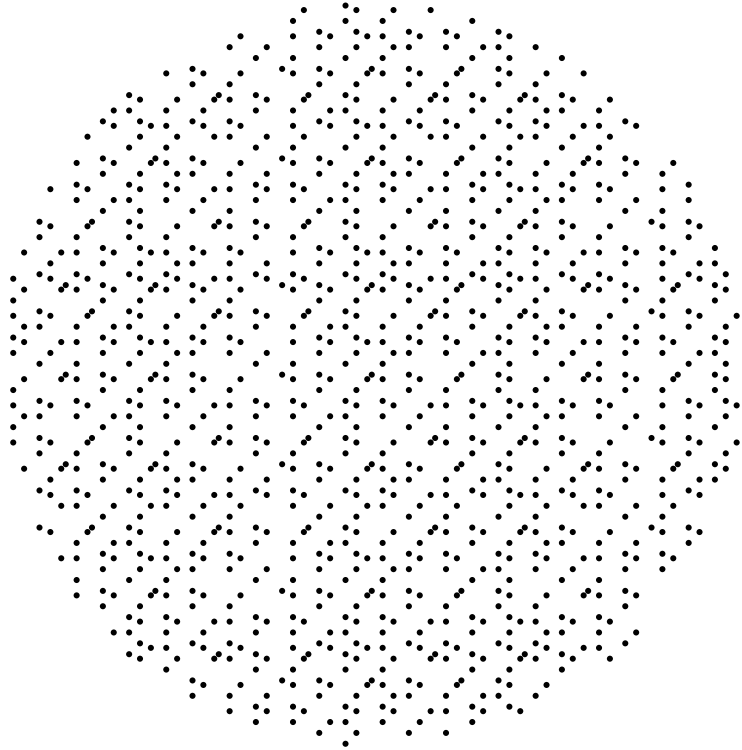}
\caption{Circular patch of the model set $\gL_{1}$.
\label{fig3}\bigskip}
\end{figure}

\begin{figure}
\includegraphics[width=0.5\textwidth]{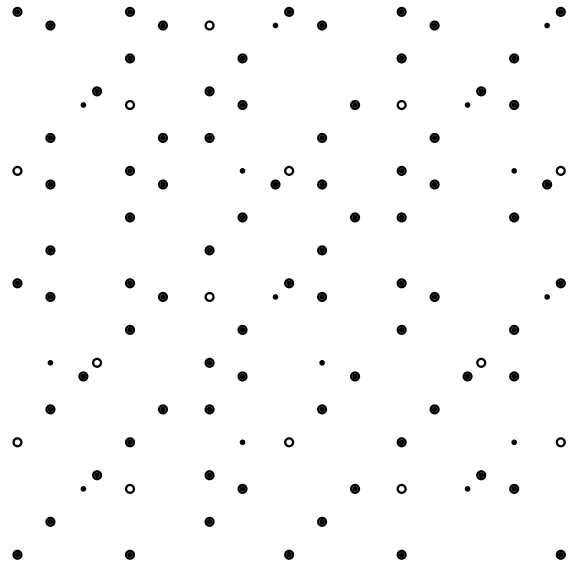}
\caption{Comparison of the two patterns $\gL_{1}$ and $\gL_{2}$. Big
dots mark points that are common to both, while open circles (small
dots) are points of $\gL_2\setminus\gL_1$ (of $\gL_1\setminus\gL_2$).
\label{fig4}}
\end{figure}

\begin{figure}
\includegraphics[width=60mm]{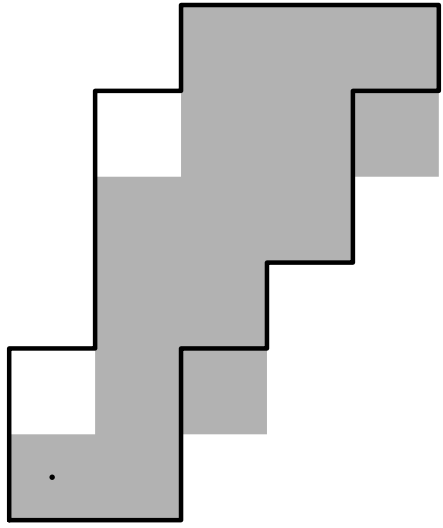}
\caption{Comparison of the two windows $P_1$ (grey area) and $P_2$
(bounded by black line).  The two white squares on the left (the two
grey ones on the right) form the window for the points in
$\gL_2\setminus\gL_1$ (in $\gL_1\setminus\gL_2$).
\label{fig5}}
\end{figure}

A patch of $\gL_1$ is shown in Figure~\ref{fig3}, while
Figure~\ref{fig4} displays the difference between the two model sets.
This can be understood in terms of the windows, compare
Figure~\ref{fig5}. In particular, the difference sets are model sets
themselves, with
\[
\dens(\gL_2\setminus\gL_1) \; =\; \dens(\gL_1\setminus\gL_2) \; =\; 
\tfrac{2}{15}\, \dens(\gL_1)\;=\; \tfrac{1}{2}.
\]
Note, however, that the sets $\gL_2\setminus\gL_1$ and
$\gL_1\setminus\gL_2$ are not homometric.

\begin{figure}
\includegraphics[width=0.62\textwidth]{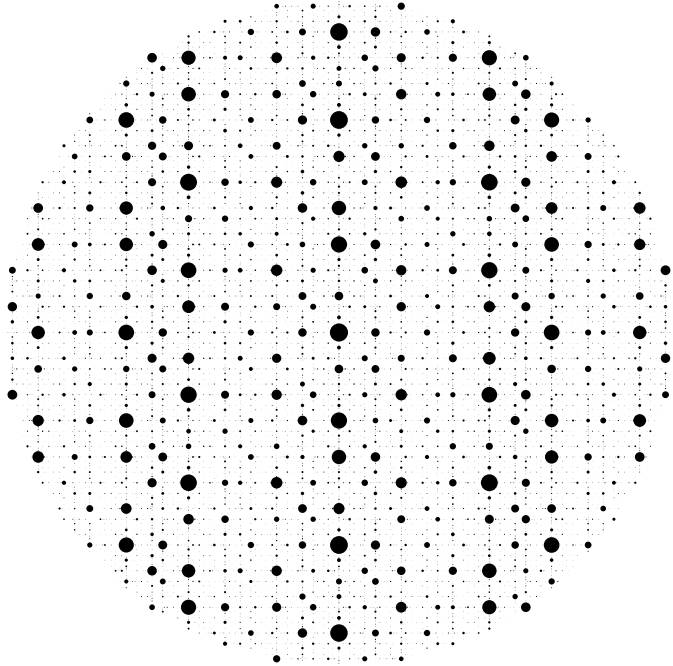}
\caption{Circular portion of the diffraction pattern, valid for both 
$\gL_1$ and $\gL_2$.  A Bragg peak at $k$ is represented by a disk 
centred at $k$ with an area proportional to the intensity.  Note 
that the image is inversion symmetric, but has no reflection 
symmetry in a line.\label{fig6}}
\end{figure}

The lattice dual to $\LL$ is denoted by $\LL^{\ast}$ (observe the
different star) and is spanned by the basis matrix
\[
B_{\LL^{\ast}}^{} \; =\; \bigl( B_{\LL}^{-1} \bigr)^{t} 
\; =\; \tfrac{1}{2}B_{\LL}^{} \, .
\]
Its projection, $\pi(\LL^{\ast})=\frac{1}{2}L$, is then the Fourier
module that carries the Bragg peaks. The diffraction measure
$\widehat{\gamma}$ is the same for both $\gL_1$ and $\gL_2$, and reads
\cite{B,H,Martin2}
\[
\widehat{\gamma}\; =\; \sum_{k\in\frac{1}{2}L} I(k^{\star})\, \delta_{k}\, .
\]
The intensity function is given by 
\[
I(x) \; =\; {\bigl(\dens(\gL_j)\bigr)}^{2}\, 
\left| \frac{1}{\vol(P_j)} \int_{P_j} e^{2\pi ixy}\; \mathrm{d} y\right|^{2}
\]
where $j\in\{1,2\}$, with $I(x)$ being the same function for both choices.
Recall
\[
\int_{a}^{b} e^{2\pi irs}\; \mathrm{d} s \; = \; 
e^{\pi i r (a+b)}\,\frac{\sin(\pi r(b-a))}{\pi r} ,
\]
the linearity of Fourier transform, and its factorisation over
rectangular domains.  With $k^{\star}=(\kappa,\lambda)$, one obtains
\[
I(\kappa,\lambda) \; = \; \frac{f(\kappa,\lambda)}{16} \; \left(
\frac{\sin(\pi\kappa)\sin(\pi\lambda)}{\pi^{2}\kappa\lambda}\right)^{2}
\]
where
\[
\begin{split}
f(\kappa,\lambda) \; = \; &
\bigl(3 + 
      2 \cos(2\pi\lambda) + 4 \cos(\pi\lambda) 
   \cos(\pi(2\kappa + 3\lambda))\bigr)\\
&\bigl(5 + 6 \cos(2\pi\kappa) + 2 \cos(4\pi\kappa) + \mbox{}\\
&\;\;  4\, (2 \cos(\pi\kappa) + \cos(3\pi\kappa)) 
  \cos(\pi (3\kappa + 6 \lambda))\bigr).
\end{split}
\]
The diffraction image is shown in Figure~\ref{fig6}.

\section{Discussion and outlook}

The characterisation of homometric model sets is possible via the
window covariogram. However, the corresponding homometry classes do
not consist of model sets only, since one can always add and subtract
point sets of density $0$. More interestingly, there is still the
possibility for other sets in these classes that differ from standard
model sets in an essential way. A better understanding of the
homometry classes defined by model sets and related structures with
pure point diffraction thus requires further investigation.

When considering homometry for systems with mixed spectra, in
particular those including diffuse scattering, the situation gets
significantly more involved. As an example from \cite{HB}, we mention
the binary Rudin-Shapiro chain in one dimension, which is
deterministic, versus the binary Bernoulli chain, which is
stochastic. They define a homometric pair of point sets that even
differ in entropy, being $0$ versus $\log(2)$ in this
example. Consequently, new concepts for further progress are required.

\section*{Acknowledgment}

It is a pleasure to thank Gabriele Bianchi and Richard Gardner for
cooperation and helpful discussions. This work was supported by the
German Research Council, within the Collaborative Research Centre 701,
and by EPSRC via Grant EP/D058465.

\end{document}